\documentclass[11pt,epsf]{article}
\usepackage{graphicx}
\usepackage{amsthm}
\usepackage{amsfonts}
\usepackage{amssymb}
\title{A property of discriminants}
\author{Vladimir Petrov Kostov\\ Universit\'e C\^ote d'Azur, CNRS, LJAD, France
\\  
e-mail: kostov@math.unice.fr} 
\date{}
\newtheorem{tm}{Theorem}

\newtheorem{rems}[tm]{Remarks}
\newtheorem{lm}[tm]{Lemma}
\newtheorem{ex}[tm]{Example}

\begin{document} 
\maketitle 
\begin{abstract}
For the family $P:=x^n+a_1x^{n-1}+\cdots +a_n$ of complex 
polynomials in the variable $x$ we study its {\em discriminant} 
$R:=$Res$(P,P',x)$, $R\in \mathbb{C}[a]$, $a=(a_1,\ldots ,a_n)$. 
When $R$ is regarded as a polynomial in $a_k$, one can consider its 
discriminant $\tilde{D}_k:=$Res$(R,\partial R/\partial a_k,a_k)$. We show that 
$\tilde{D}_k=c_k(a_n)^{d(n,k)}M_k^2T_k^3$, where $c_k\in \mathbb{Q}^*$, 
$d(n,k):=\min (1,n-k)+\max (0,n-k-2)$, the polynomials 
$M_k,T_k\in \mathbb{C}[a^k]$ have integer coefficients, 
$a^k=(a_1,\ldots ,a_{k-1},a_{k+1},\ldots ,a_n)$, the sets $\{ M_k=0\}$ 
and $\{ T_k=0\}$ are the projections in the space of the variables $a^k$ of the 
closures of the strata of the variety 
$\{ R=0\}$ on which $P$ has respectively two double 
roots or a triple root. Set $P_k:=P-xP'/(n-k)$ for $1\leq k\leq n-1$ and 
$P_n:=P'$. One has $T_k=${\rm Res}$(P_k,P_k',x)$ for $k\neq n-1$ and 
$T_{n-1}=${\rm Res}$(P_{n-1},P_{n-1}',x)/a_n$.\\ 

{\bf AMS classification:} 12E05; 12D05\\

{\bf Key words:} polynomial in one variable; discriminant set; 
resultant; multiple root
\end{abstract}  

\section{Introduction}

In the present paper we consider the general family of monic degree $n$ 
complex polynomials in one variable $P:=x^n+a_1x^{n-1}+\cdots +a_n$. (For 
$a_1=0$ this is the versal deformation of the $A_{n-1}$-singularity, 
see~\cite{AGV}). Its 
{\em discriminant} is the resultant $R:=$Res$(P,P',x)$, i.e. the determinant of 
the {\em Sylvester matrix} $S(P,P',x)$. We remind that $S(P,P',x)$ is 
$(2n-1)\times (2n-1)$, its first (resp. $n$th) row equals 

$$(1,a_1,\ldots ,a_n,0,\ldots ,0)~~~~\, \, \, \, {\rm (resp.}~~\, 
(n,(n-1)a_1,\ldots ,a_{n-1},0,\ldots ,0)~~{\rm )}~,$$
its second (resp. $(n+1)$st) row is obtained by shifting the first 
(resp. the $n$th) one to the right by one position while adding $0$ 
to the left etc. 
Set 
$a:=(a_1,\ldots ,a_n)$, $a^k:=(a_1,\ldots,a_{k-1},a_{k+1},\ldots ,a_n)$ and 
$R_{a_k}:=\partial R/\partial a_k$. It is well-known that: 
\vspace{1mm}

A) $R$ is a quasi-homogeneous polynomial in the coefficients $a_j$, 
where the quasi-homogeneous weight of $a_j$ equals $j$. It is a degree $n$ 
polynomial in each of the variables $a_j$, $1\leq j\leq n-1$, and a degree 
$n-1$ polynomial in $a_n$.
\vspace{1mm}

B) The set $\{ R=0\}$ is the set of values of the coefficients $a_j$ for which 
$P$ has a multiple root. It contains the subsets $\Sigma$ and $\tilde{M}$ 
(the {\em Maxwell stratum}) such that for 
$a\in \Sigma$ (resp. $a\in \tilde{M}$) the polynomial $P$ 
has a root of multiplicity $3$ (resp. has two different double roots). 
The semi-algebraic sets $\Sigma$ and $\tilde{M}$ are irreducible. 
Indeed, the closure of 
$\Sigma$ is the 
image of the map $\mathbb{C}^{n-2}\rightarrow \mathbb{C}^{n-2}$, 
$(z_1,z_4,z_5,\ldots ,z_n)\mapsto a$, where in the computation of 
$(-1)^ja_j$ as $j$th elementary symmetric function of $z_1$, $\ldots$, $z_n$ 
one sets 
$z_2=z_3=z_1$; the closure of $\tilde{M}$ is the 
image of the map $\mathbb{C}^{n-2}\rightarrow \mathbb{C}^{n-2}$, 
$(z_1,z_3,z_5,z_6\ldots ,z_n)\mapsto a$, 
where in the computation of $a$ one sets 
$z_2=z_1$ and $z_4=z_3$. It is easy to see that the intersections of the sets 
$\Sigma$ and $\tilde{M}$ with each of the subspaces $\{ a_j=0\}$ are 
proper subsets of $\Sigma$ and $\tilde{M}$.  
\vspace{1mm}

One can consider $R$ as a polynomial in $a_k$, with 
coefficients in $\mathbb{C}[a^k]$.  
Thus one is led to consider the repeated resultants 
$\tilde{D}_k:=$Res$(R,R_{a_k},a_k)$. The following result 
is proved in~\cite{Ko1} (see Proposition~7 there):

\begin{lm}\label{lmdivisible}
Set $d(n,k):=\min (1,n-k)+\max (0,n-k-2)$. The polynomial 
$\tilde{D}_k$ equals $(a_n)^{d(n,k)}\tilde{D^0_k}$, where 
$\tilde{D^0_k}\in \mathbb{C}[a]$ is not divisible by any of the 
variables $a_i$, $1\leq i\leq n$.
\end{lm}

%
%
%

\begin{ex}\label{exn=3}
{\rm For $n=3$ one has $P:=x^3+ax^2+bx+c$, $P'=3x^2+2ax+b$ and 

$$R:={\rm Res}(P,P',x)=4a^3c-a^2b^2-18abc+4b^3+27c^2~.$$ 
Set $\tilde{D}_a:=$Res$(R,\partial R/\partial a,a)$ 
and similarly for $\tilde{D}_b$ and $\tilde{D}_c$. 
Hence} 

$$\tilde{D}_a=-64c(b^3-27c^2)^3~~,~~
\tilde{D}_b=-64c(a^3-27c)^3~~{\rm and}~~
\tilde{D}_c=-432(-3b+a^2)^3~.$$
\end{ex}

\begin{ex}\label{exn=4}
{\rm For $n=4$ one has $P:=x^4+ax^3+bx^2+cx+d$, $P'=4x^3+3ax^2+2bx+c$ and} 

$$\begin{array}{ccccl}
R&:=&{\rm Res}(P,P',x)&=&
-27a^4d^2+18a^3bcd-4a^3c^3+a^2b^2c^2+144a^2bd^2-4a^2b^3d\\
&&&&-6a^2c^2d-80ab^2cd+18abc^3-192acd^2+16b^4d\\ 
&&&&-4b^3c^2-128b^2d^2+144bc^2d-27c^4+256d^3~.\end{array}$$
{\rm One finds that} 

$$\tilde{D}_a=6912d^2M_a^2T_a^3~~,~~\tilde{D}_b=-4096dM_b^2T_b^3~~,~~
\tilde{D}_c=6912dM_c^2T_c^3~~{\rm and}~~\tilde{D}_d=4096M_d^2T_d^3~~,$$
{\rm where the factors $M_a$, $T_a$, $M_b$, $\ldots$, $T_d$ are irreducible:}

$$\begin{array}{lcl}
M_a=16b^2d^2-8bc^2d+c^4-64d^3&,&T_a=
3b^4d-b^3c^2+72b^2d^2-108bc^2d+27c^4+432d^3\\ \\ 
M_b=a^2d-c^2&,&T_b=27a^4d^2-a^3c^3-6a^2c^2d-768acd^2+27c^4+4096d^3\\ \\ 
M_c=a^4-8a^2b+16b^2-64d&,&T_c=
27a^4d-a^2b^3-108a^2bd+3b^4+72b^2d+432d^2\\ \\ 
M_d=a^3-4ab+8c&,&T_d=27a^3c-9a^2b^2-108abc+32b^3+108c^2~.\end{array}$$
{\rm One can notice that the equation $M_b=0$ defines the Whitney umbrella.}
\end{ex}

We prove the following theorem:

\begin{tm}\label{maintmreduc}
For $n\geq 4$ the polynomial $\tilde{D}_k$ is of the form 
$c_k(a_n)^{d(n,k)}M_k^2T_k^3$, where $c_k\in \mathbb{Q}^*$, 
the degree $d(n,k)$ is defined in  
Lemma~\ref{lmdivisible} 
and the polynomials $M_k, T_k\in \mathbb{C}[a^k]$ are with integer coefficients 
and irreducible. The zero sets of these polynomials are the closures of the 
projections in the space of the variables $a^k$ of the sets $\tilde{M}$ and 
$\Sigma$. 
\end{tm}

The proofs of Theorem~\ref{maintmreduc}, Lemma~\ref{nn+1} and Lemma~\ref{QHD} 
are to be found in Section~\ref{secproofs}.

{\bf Acknowledgement.} The author is grateful to B.Z. Shapiro from the 
University of Stockholm for the formulation of the problem and its 
subsequent discussions.

\section{Comments and lemmas}

Theorem~\ref{maintmreduc} is formulated for $n\geq 4$ because for $n<4$ the set 
$\tilde{M}$ does not exist. In Example~\ref{exn=3} only the cubes of the 
factors $T_k$ and the powers of $a_n$ (i.e. of $c$) are present.

It is well-known that 
$R=\prod_{1\leq i<j\leq n}(z_i-z_j)^2$. Denote by $\Delta$ the union of 
hyperplanes $\{ z_i=z_j\}$ in the space $\mathbb{C}^n$ of the roots of 
the polynomial $P$. In the last presentation of $R$ as a product 
it is necessary to have the differences of roots $z_i-z_j$ squared 
because when the roots change continuously along a loop avoiding the set 
$\Delta$ so that in the end two of them are exchanged, then such an exchange 
should not change the value of $R$.  

By analogy, the fact that the power of the factor $T_k$ in the formula 
for $\tilde{D}_k$ in Theorem~\ref{maintmreduc}) is a multiple of  
$3$ can be explained like this. 
At a point $a=a^*\in \Sigma$ (we assume that 
$a^*\not\in \bar{\Sigma}\backslash \Sigma$) 
three roots 
$z_1$, $z_2$, $z_3$ of $P$ coalesce. For fixed nearby values of 
$a^k$ the polynomial $R$ (when considered as a polynomial in $a_k$) 
has two roots $\zeta _1$ and $\zeta _2$ that coalesce for $a^k={a^*}^k$ 
(the projection of $a^*$ in the space of the variables $a^k$). 
These roots correspond  to equalities and inequalities between the roots of $P$ 
of the form $z_1=z_2\neq z_3$ and $z_1\neq z_2=z_3$ for $a^k\neq {a^*}^k$, 
and to $z_1=z_2=z_3$ for $a^k={a^*}^k$. When the 
$(n-1)$-tuple of coefficients $a^k$ circumvents 
the projection $\Sigma _k$ of $\Sigma$ in the space of the variables $a^k$ 
along a generic loop, the 
three roots $z_i$ of $P$ undergo a cyclic permutation of order $3$ 
and now the roots 
$\zeta _1$ and $\zeta _2$ of $R$ correspond to other equalities and 
inequalities between the roots $z_i$, namely, to $z_3=z_1\neq z_2$ 
and $z_3\neq z_1=z_2$. In order $\tilde{D}_k$ 
to be invariant w.r.t. such 
permutations the power of $T_k$ dividing the resultant $\tilde{D}_k$ 
must be a multiple of $3$.

\begin{figure}[htbp]
\centerline{\hbox{\includegraphics[scale=0.7]{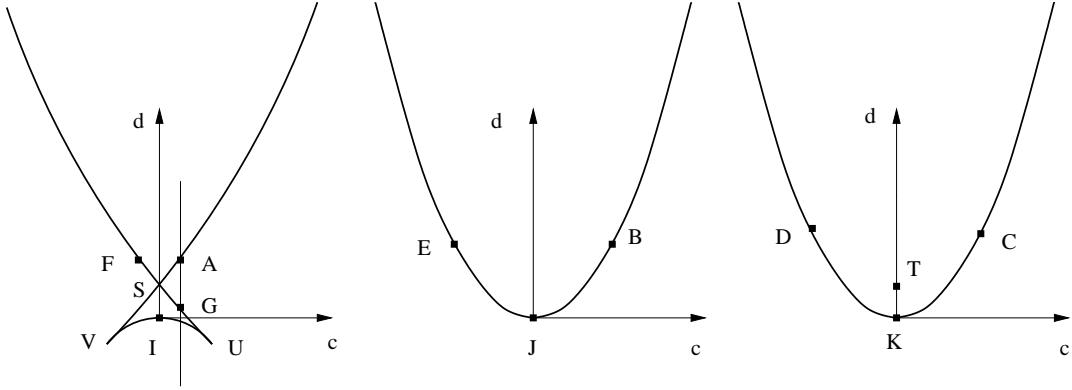}}}
    \caption{The sets $\{ R=0\} |_{a=0,b=-1}$, $\{ R=0\} |_{a=b=0}$ 
and $\{ R=0\} |_{a=0,b=1}$ for $n=4$.}
\label{ws}
\end{figure}

For the power of $M_k$ being even a similar explanation exists. To this 
end we remind first some facts about $R$ for $n=4$. The formula for $R$ 
was obtained in Example~\ref{exn=4}. On Fig.~\ref{ws} we show for real values 
of $c$ and $d$ 
the sets $\{ R=0\} |_{a=0,b=-1}$, $\{ R=0\} |_{a=b=0}$ 
and $\{ R=0\} |_{a=0,b=1}$ (from left to right) which are symmetric w.r.t. the 
$d$-axis. This figure can be compared with the well-known picture 
of the swallowtail catastrophe, see \cite{PS}.
Fig.~\ref{ws} gives a sufficient idea about the set 
$\{ R=0\} |_{a=0}$ because the set $\{ R=0\}$ is invariant under the 
quasi-homogeneous dilatations $a\mapsto ta$, $b\mapsto t^2b$, $c\mapsto t^3c$, 
$d\mapsto t^4d$, $t\neq 0$. 

At the points $U$ and $V$ the polynomial $P$ 
has a triple real and a simple real root ($U$ and $V$ are ordinary $2/3$-cusp 
points for the real curve $\{ R=0\} |_{a=0,b=-1}$). One has 

$$\Sigma \cap \{ a=0,b=-1\} =\{ U,V\}~~~,~~~
\tilde{M}\cap \{ a=0,b=-1\} =\{ S\} ~.$$ 
At the point $S$ (with $d$-coordinate equal to $1/4$) 
the curve $\{ R=0\} |_{a=0,b=-1}$ has transversal self-intersection  
and the polynomial $P$ has two double real roots. At the point $T$ (which is an 
isolated double point of the real curve $\{ R=0\} |_{a=0,b=1}$, with 
$d$-coordinate equal to $1/4$) the polynomial 
$P$ has a double complex conjugate pair. At the points $I$, $J$ and $K$ one has 
$c=d=0$. The real curves $\{ R=0\} |_{a=0,b=-1}$ 
and $\{ R=0\} |_{a=0,b=1}$ are smooth at $I$ and $K$ respectively while 
$\{ R=0\} |_{a=b=0}$ has a $4/3$-type singularity at $J$.

From now on we keep in mind that the set $\{ R=0\}$ can be defined in both 
contexts -- the ones of real or of complex variables $x$, $a$, $b$, $c$ and 
$d$. In this sense we make use of Fig.~\ref{ws} as an illustration of the 
real case and as a hint for the complex one. Why for $n=4$ the powers of 
the factors $M_k$ should be even is suggested by the following lemma. For 
$n>4$ the analogs of the loops $\bar{\gamma}$ and $\Gamma$ of the lemma 
exist in a 
neighbourhood of any value of the parameters $a_j$ for which the polynomial $P$ 
has a quadruple root, but their explicit construction is harder to describe. 

\begin{lm}\label{lmpath}
In the complex case there exists a loop $\bar{\gamma}$ 
belonging to the space of 
variables $(b,c)$ which can be lifted to a loop  
$\Gamma \subset \{ R=0\} |_{a=0}$ circumventing the set 
$\overline{\Sigma \cup \tilde{M}}$ such that any fibre of the projection 
$\Gamma \rightarrow \bar{\gamma}$ consists of two points and 
the monodromy defined on the fibre after one turn along $\bar{\gamma}$ 
is nontrivial.
\end{lm}

\begin{proof}
In what follows an additional index $d$ denotes the projection of a given 
set in the space of variables $(b,c,d)$ ($a$ is presumed equal to $0$) into 
the space of variables $(b,c)$. 
Consider the point $A$ on Fig.~\ref{ws}. We are going 
to construct a continuous path $\gamma \subset \{ R=0\} |_{a=0}$ 
leading from $A$ to $G$, one of the two points of $\{ R=0\} |_{a=0}$ 
which share with $A$ the same 
$b$- and $c$-coordinates as shown on Fig.~\ref{ws}. 
As $b$ increases from $-1$ to $1$, 
the point $A$ becomes the point $B$ for $b=0$ and then $C$ for $b=1$. 
Then we decrease $c$ 
by keeping the same value of $b$ -- this gives the arc $CKD$. Then we 
fix $c$ and decrese $b$ -- this gives the arc $DEF$. Finally 
we add the arc $FG$. The thus constructed path is real. 
Three remarks will be needed for what follows:
\vspace{1mm}

1) The path $\gamma$, in its part between the points $A$ and $F$, 
can be constructed as symmetric w.r.t. the plane $\{ c=0\}$.
\vspace{1mm}

2) The projection $\Sigma _d$ of $\Sigma$ 
is defined by $32b^3+108c^2=0$, i.e. $8b^3+27c^2=0$; 
the equation of this semi-cubic 
parabola is obtained from the equation 
$T_d=0$ by setting $a=0$, see Example~\ref{exn=4}. There exists  
a unique number $b_0\in (-1,0)$ such that for $b=b_0$ the projection 
$\gamma _d$ of $\gamma$ intersects $\Sigma _d$ at two points 
$(b_0,\pm c_0)$. 
\vspace{1mm}

3) In the real case the path $\gamma$ has to pass through the point 
$S\in \tilde{M}$, 
but in the complex one $\gamma$ can be modified so that it circumvent $S$. 
The points of 
the modified path $\gamma$ which are close to $S$ do not have 
all their coordinates real.   
\vspace{1mm}

Now we construct (in the complex case) a path 
$\gamma ^1 \subset \{ R=0\} |_{a=0}$ leading from $G$ to $A$ and satisfying the 
condition $\gamma ^1_d=\gamma _d$. At the same time 
we modify the path $\gamma$ in order to have this condition. 
If the path $\gamma _1$ is defined such that $\gamma ^1_d=\gamma _d$, 
then for $b=b_0$, $\gamma _1$ will intersect 
the set $\Sigma$. Therefore for $b$ close to $b_0$ we modify $\gamma _1$ 
and $\gamma$ so that $\gamma ^1$ avoid the set $\Sigma$. 
(We make two such modifications, 
corresponding to points of $\gamma _d$ and $\gamma ^1_d$ close to 
$(b_0,c_0)$ and to $(b_0,-c_0)$. The modifications can be made 
symmetrically w.r.t. the plane $\{ c=0\}$.)

For the values of $b$ close to $b_0$ the points 
of $\gamma$ do not have all their coordinates real. As for $\gamma ^1$, its 
points do not have all coordinates real not only for $b$ close to $b_0$, but 
also for $b\in [b_0,1]$ (recall the construction of the 
arcs $ABC$ and $DEF$ of $\gamma$) and for 
$b=1$, $c\neq 0$ (recall the construction of its arc $CKD$). 
Indeed, as $R$ is a degree 
$3$ polynomial in $d$, then in the real case  
it has either three real roots (see for instance 
the vertical 
line on the left part of Fig.~\ref{ws} which intersects the set $\{ R=0\}$ 
at three points two of which are $A$ and $G$) or one real 
and two complex conjugate ones; this is, 
in particular, the case of any vertical line different from the $d$-axis for 
$b=1$, see the right part of Fig.~\ref{ws}. (The $d$-axis on the right part of 
the figure 
corresponds to one simple root at $0$ and a double one at $1/4$. One simple 
and one double real root is also the situation observed on the vertical lines 
passing through the points $U$ and $V$.) 

To obtain the proof of the lemma one sets $\bar{\gamma}=\gamma _d=\gamma ^1_d$ 
and one defines the loop 
$\Gamma$ as the concatenation of $\gamma$ and $\gamma ^1$. For points of 
$\gamma$ and $\gamma ^1$ close to the point $S$ one has $\gamma _d=\gamma ^1_d$ 
and no self-intersection of $\Gamma$ takes place.
\end{proof}

\begin{rems}\label{remsimportant}
{\rm (1) To prove Theorem~\ref{maintmreduc} we need to recall 
some notation and results from \cite{Ko1}. 
Suppose that $G_1$ and $G_2$ are polynomials in several variables 
one of which is denoted by $y$. By $S(G_1,G_2,y)$ we denote 
the Sylvester matrix of $G_1$ and $G_2$ when considered as polynomials in $y$. 
We set $P_k:=P-xP'/(n-k)$ for $1\leq k\leq n-1$ and $P_n:=P'$.

(2) It is shown in \cite{Ko1} that for $k\neq n-1$ 
the polynomial $V_k:=${\rm Res}$(P_k,P_k',x)$ is irreducible and that the 
polynomial {\rm Res}$(P_{n-1},P_{n-1}',x)$ is the product of $a_n$ 
and an irreducible polynomial in~$a^{n-1}$. We set 
$V_{n-1}:=${\rm Res}$(P_{n-1},P_{n-1}',x)/a_n$. It follows from Theorem~12 
of \cite{Ko1} that $V_k=T_k$, $k=1,\ldots ,n$. 
Theorem~\ref{maintmreduc} allows to find the polynomials $M_k$ and $T_k$; 
however the definition of $T_k$ as $T_k=V_k$ is an easier way to find $T_k$.

(3) We denote by QHD$(U)$ the quasi-homogeneous degree of a quasi-homogeneous 
polynomial $U\in \mathbb{C}[a]$, where the quasi-homogeneous weight of 
$a_k$ is $k$.

(4) Set $Q_k:=(n-k)P_k=(n-k)P-xP'$, $k\leq n-1$, $Q_n:=P'$. 
When we compare polynomials $P_k$, $Q_k$, $R$ or $V_k$ 
for two consecutive values 
of $n$ (i.e. for $n$ and $n+1$) 
we write $P_k^n$, $P_k^{n+1}$, $Q_k^n$, $Q_k^{n+1}$, $R^n$, $R^{n+1}$ or 
$V_k^n$, $V_k^{n+1}$. 
Notice that as $Q_k=-kx^n+\sum _{j=1}^n(j-k)a_jx^{n-j}$, one has} 

\begin{equation}\label{QP}
Q_k^{n+1}=xQ_k^n+(n+1-k)a_{n+1}~~{\rm and}~~(Q_k^{n+1})'=x(Q_k^n)'+Q_k^n~.
\end{equation}
\end{rems}
 
In the following lemma and 
its proof $\Omega$ denotes nonspecified nonzero rational numbers. 

\begin{lm}\label{nn+1}
(1) One has $V_*:=V_k^{n+1}|_{a_{n+1}=0}=\Omega (a_n)^2V_k^n$ for $1\leq k\leq n-2$, 
$V_*=\Omega (a_n)^3V_k^n$ for $k=n-1$ and 
$V_*=\Omega (a_{n-1})^3V_k^n$ for $k=n$.

(2) One has $R^{n+1}|_{a_{n+1}=0}=\pm a_n^2R^n$. 
\end{lm}

The following lemma announces the quasi-homogeneous degrees of certain 
polynomials that appear in this text:

\begin{lm}\label{QHD}
For $n\geq 4$ one has the following quasi-homogeneous degrees of polynomials:

(1) {\rm QHD}$(R)=${\rm QHD}$(V_k)=n(n-1)$, $1\leq k\leq n-2$.
 
(2) {\rm QHD}$(V_{n-1})=n(n-2)$. 

(3) {\rm QHD}$(V_n)=(n-1)(n-2)$.

(4) {\rm QHD}$(R_{a_k})=n(n-1)-k$, $1\leq k\leq n-2$, 
{\rm QHD}$(R_{a_{n-1}})=n^2-3n+1$, {\rm QHD}$(R_{a_n})=n^2-4n+2$. 

(5) {\rm QHD}$(\tilde{D}_k)=n(n-1)^2+n^2(n-k-1)$, $1\leq k\leq n-1$, 
{\rm QHD}$(\tilde{D}_n)=n(n-1)(n-2)$.

(6) {\rm QHD}$(M_k)=n^3-3n^2+2n-(n^2-n)(k+1)/2$, $1\leq k\leq n-2$, 
{\rm QHD}$(M_{n-1})=n(n-2)(n-3)/2$, {\rm QHD}$(M_n)=(n-1)(n-2)(n-3)/2$.
\end{lm}

\section{Proofs\protect\label{secproofs}}

\begin{proof}[Proof of Lemma~\ref{nn+1}] 
The equality $A=[B]_{\ell ,r}$ means that the matrix $A$ is obtained from the 
matrix $B$ by deleting its $\ell$th row and $r$th column. 
Prove part (1). In the proof of the lemma we use the polynomials $Q_k$ 
instead of $P_k$. 
For $1\leq k\leq n-2$ set $Q_*:=Q_k^{n+1}|_{a_{n+1}=0}=xQ_k^n$. 
Consider the $(2n+1,2n+1)$-Sylvester matrix 
$S_*:=S(Q_*,Q_*',x)$. The only nonzero entry in its last column is 
$\Omega a_n$ in position $(2n+1,2n+1)$. Hence when finding its determinant 
$\Omega V_*$ one can 
develop it w.r.t. the last column to obtain $V_*=\Omega a_nV_{**}$, 
where $V_{**}=\det S_{**}$, $S_{**}=[S_*]_{2n+1,2n+1}$.

Subtract for $j=1,\ldots ,n$ the $j$th row of $S_{**}$ 
from its $(n+j)$th row. 
This doesn't change $V_{**}$.  
Hence the terms $\Omega a_n$ disappear in the 
$(n+1)$st, $\ldots$, $(2n)$th rows 
of $S_{**}$, see (\ref{QP}). 
The only nonzero entry of the 
new matrix (denoted by $S_{***}$) 
in its last column is $\Omega a_n$ in position $(n,2n)$. It is easy to see that 
$[S_{***}]_{n,2n}=S(Q_k^n,(Q_k^n)',x)$ (this can be deduced from (\ref{QP})). 
Hence $V_{**}=\det S_{***}=\Omega a_nV_k^n$ and $V_*=\Omega (a_n)^2V_k^n$.

For $k=n-1$ the above reasoning differs only in the end -- one defines 
$V_{n-1}^n$ not as $\det ([S_{***}]_{n,2n})$ 
(the latter is divisible by $a_n$), but as 
$\det ([S_{***}]_{n,2n})/a_n$. Hence $V_*=\Omega (a_n)^3V_{n-1}^n$.

For $k=n$ consider the $(2n+1)\times (2n+1)$-matrix 
$S^0:=S(Q_n^{n+1},(Q_n^{n+1})',x)$. Its last column contains 
a single nonzero entry 
($\Omega a_{n+1}$ in position $(n,2n+1)$). 
By definition $\det S^0=\Omega a_{n+1}V_n^{n+1}$. 
Hence $V_*=\Omega \det S^{\dagger}$, where 
$S^{\dagger}=([S^0]_{n,2n+1})|_{a_{n+1}=0}$. 

The last column of $S^{\dagger}$ contains a single nonzero entry 
($\Omega a_{n-1}$ in 
position $(2n,2n)$), so to find $\det S^{\dagger}$ one can develop it w.r.t. 
the last column. This gives $V_*=\Omega a_{n-1}\det S^{\dagger 0}$, where 
$S^{\dagger 0}=[S^{\dagger}]_{2n,2n}$. 

Subtract the $j$th row of $S^{\dagger 0}$ from its $(n-1+j)$th one, 
$j=1,\ldots ,n-1$; hence the terms $\Omega a_{n-1}$ disappear in the 
$n$th, $\ldots$, $(2n-2)$nd rows (see (\ref{QP})).
This gives the matrix $S^{\dagger *}$ such that 
$\det S^{\dagger *}=\det S^{\dagger 0}$. 

The only nonzero entry in the last column of $S^{\dagger *}$ is 
$\Omega a_{n-1}$ in position $(2n-1,2n-1)$. Hence 
$\det S^{\dagger *}=\Omega a_{n-1}\det S^{\dagger \dagger}$, where 
$S^{\dagger \dagger}=[S^{\dagger *}]_{2n-1,2n-1}$. The only nonzero entry 
of $S^{\dagger \dagger}$ in its last column 
is in position $(n-1,2n-2)$ and equals $\Omega a_{n-1}$. Thus 
$V_*=\Omega (a_{n-1})^3\det S^{\dagger \dagger 0}$, where 
$S^{\dagger \dagger 0}=[S^{\dagger \dagger}]_{n-1,2n-2}$. The $(2n-3)\times (2n-3)$-matrix 
$S^{\dagger \dagger 0}$ equals $S(Q_n^n/x,(Q_n^n/x)',x)$, i.e. 
$\Omega S((P^n)',(P^n)'',x)$.  

To prove part (2) one notices that for $a_{n+1}=0$ one has $P^{n+1}=xP^n$ 
and the Sylvester matrix 
$S^1:=S(xP^n,(xP^n)',x)$ contains a single nonzero entry in its last column, 
namely $a_n$ in position $(2n+1,2n+1)$. Set $S^2:=[S^1]_{2n+1,2n+1}$. Hence 
$R^{n+1}|_{a_{n+1}=0}=\det S^1=a_n\det S^2$. For $j=1$, $\ldots$, $n$ subtract the 
$j$th row of $S^2$ from its $(n+j)$th one. The newly obtained matrix 
(denoted by $S^3$) has a single nonzero entry in its last column. This is 
$a_n$ in position $(n,2n)$. Set $S^3:=[S^2]_{n,2n}$. Hence 
$\det S^2=\pm a_n\det S^3$, i.e. $R^{n+1}|_{a_{n+1}=0}=\pm a_n^2\det S^3$. 
On the other hand $S^3=S(P^n,(P^n)',x)$ from which part (2) follows.     

\end{proof}

\begin{proof}[Proof of Lemma~\ref{QHD}]
We denote by $W$ any of the polynomials $R$, $V_k$, $k\leq n-2$, or 
$a_nV_{n-1}$ and we remind that $T_k=V_k$, see Remarks~\ref{remsimportant}. 
Any polynomial $W$ 
contains a monomial $\beta a_n^{n-1}$, $\beta \neq 0$. Indeed, the only 
positions in which the matrix $S(W,W',x)$ contains the variable $a_n$ 
are $(i,n+i)$, $i=1,\ldots ,n-1$; in these 
positions the matrix has terms of the form $\eta a_n$, $\eta \neq 0$. When 
$\det (S(W,W',x))$ is computed, these 
terms are multiplied by the constant nonzero terms in positions 
$(n-1+j,j)$, $j=1,\ldots ,n$ to give the only monomial of the form 
$\beta a_n^{n-1}$ in $\det (S(W,W',x))$. 
Hence QHD$(R)=$QHD$(V_k)=$QHD$(a_nV_{n-1})=n(n-1)$ which proves parts (1) and 
(2). The proof of part (3) is analogous (one considers polynomials $W$ of 
degree $n-1$ instead of $n$ and $a_{n-1}$ plays the role of $a_n$).  

Part (4) follows from parts (1), (2) and (3) -- when $R$ is 
differentiated w.r.t. $a_k$, its quasi-homogeneous degree decreases by $k$. 

Prove part (5). For $a_i=0$, $k\neq i\neq n$, $k<n$, one has 
$R=\Omega _1a_k^na_n^{n-k-1}+\Omega _2a_n^{n-1}$, 
$\Omega _1\neq 0\neq \Omega _2$, see Statement~8 in~\cite{Ko1}. Therefore 
the Sylvester matrix $S(R,R_{a_k},a_k)$ has only the following nonzero entries, 
in the following positions:

$$\begin{array}{clcclcl}
\Omega _1a_n^{n-k-1}&{\rm at}~~(i,i)&,&
\Omega _2a_n^{n-1}&{\rm at}~~(i,n+i)&,&i=1,\ldots ,n-1\\ 
&&{\rm and}&n\Omega _1a_n^{n-k-1}&{\rm at}~~(n-1+j,j)&,&j=1,\ldots ,n~.
\end{array}$$
Hence its determinant equals $\Omega a_n^{(n-1)^2+n(n-k-1)}$, $\Omega \neq 0$ 
which proves part (5) for $k<n$. 

If $k=n$ and $a_i=0$ for $i\leq n-2$, then 
$R=\Omega _3a_n^{n-1}+\Omega _4a_{n-1}^n$, $\Omega _3\neq 0\neq \Omega _4$. 
Indeed, the presence of the monomials $\Omega _3a_n^{n-1}$ and 
$\Omega _4a_{n-1}^n$ in $R$ is easy to deduce from the form of the matrix 
$S(P,P',x)$, and for $a_i=0$ ($i\leq n-2$) there exist no other monomials 
of quasi-homogeneous weight $n(n-1)$ in Res$(P,P',x)$. Hence 
the Sylvester matrix $S(R,R_{a_n},a_n)$ (of size $(2n-3)\times (2n-3)$) 
has only the following nonzero entries, 
in the following positions:

$$\begin{array}{clcclcl}
\Omega _3&{\rm at}~~(i,i)&,&
\Omega _4a_{n-1}^n&{\rm at}~~(i,n-1+i)&,&i=1,\ldots ,n-2\\ 
&&{\rm and}&(n-1)\Omega _3&{\rm at}~~(n-2+j,j)&,&j=1,\ldots ,n-1~.
\end{array}$$
Hence its determinant equals 
$\tilde{\Omega}a_{n-1}^{n(n-2)}$, $\tilde{\Omega}\neq 0$. Part (5) is proved. 

Part (6) follows from the previous parts, from Lemma~\ref{lmdivisible} and 
from Theorem~\ref{maintmreduc}. Indeed, for $k\leq n-2$ one has 

$$\begin{array}{ccl}
{\rm QHD}(M_k)&=&({\rm QHD}(\tilde{D}_k)-3{\rm QHD}(V_k)-n(n-k-1))/2\\
&=&(n(n-1)^2+n^2(n-k-1)-3n(n-1)-n(n-k-1))/2\\
&=&n^3-3n^2+2n-(n^2-n)(k+1)/2~.
\end{array}$$
For $k=n-1$ one obtains 

$$\begin{array}{cclcl}
{\rm QHD}(M_{n-1})&=&({\rm QHD}(\tilde{D}_{n-1})-3{\rm QHD}(V_{n-1})-n)/2&&\\
&=&(n(n-1)^2-3n(n-2)-n)/2&=&n(n-2)(n-3)/2~.
\end{array}$$
Finally for $k=n$ one gets 

$${\rm QHD}(M_n)=({\rm QHD}(\tilde{D}_n)-3{\rm QHD}(V_n))/2=
(n-1)(n-2)(n-3)/2~.$$
    
\end{proof}

\begin{proof}[Proof of Theorem~\ref{maintmreduc}]
At a point of the set $\{ R=0\}$, where $P$ has one double nonzero root 
and $n-2$ simple roots, this set is locally the graph of a function 
analytic in the variables $a^k$, for any $1\leq k \leq n$; if the double 
root is at $0$, then this property holds for $k=n$ and fails for 
$1\leq k\leq n-1$; at a point of this set for which $P$ has a root 
of multiplicity $\geq 3$ the set is not smooth 
(see Theorem~4 in \cite{Ko1}). It is not smooth also at points 
for which $P$ has $m\geq 2$ double roots and $n-2m$ simple ones; 
at such points the set 
$\{ R=0\}$ is locally the transversal intersection of $m$ 
smooth hypersurfaces (see part (1) of Remarks~6 in \cite{Ko1}). 

Hence a priori the polynomial $\tilde{D}_k$ is of the form 
$(a_n)^{s_k}M_k^{\alpha _k}T_k^{\beta _k}$, where $s_k\in \mathbb{N}\cup 0$, 
$\alpha _k,\beta _k\in \mathbb{N}$, $\{ M_k=0\}$ (resp. $\{ T_k=0\}$) 
is the projection of the set $\tilde{M}$ (resp. of $\Sigma$) 
in the space of the variables $a^k$.
The equality $s_k=d(n,k)$ follows from Lemma~\ref{lmdivisible}. 

Further we prove the theorem by induction on $n$. For $n=4$ its proof 
follows from Example~\ref{exn=4}. Suppose that for some 
$a\in \mathbb{C}^{n+1}$ the polynomial $P^{n+1}$ has a simple root 
$h\in \mathbb{C}$. Set $x\mapsto x+h$. The new polynomial $P^{n+1}$ 
has a simple root at $0$ hence $a_{n+1}=0$. The discriminant $R^{n+1}$ 
depends only 
on the differences between the roots of $P^{n+1}$ hence it remains invariant 
under shifts of the variable $x$. For $a_{n+1}=0$ 
one can apply Lemma~\ref{nn+1}. The lemma implies that for $k\leq n-1$ the 
discriminant 
Res$(R^{n+1},\partial R^{n+1}/\partial a_k,a_k)$ is of the form 
$a_n^{t_k}M_k^2T_k^3$, $t_k\in \mathbb{N}$, i.e. one has $\alpha _k=2$ and 
$\beta _k=3$ for 
$k\leq n-1$, $a_n\neq 0$ and $a_{n-1}\neq 0$. The sets $\tilde{M}$ and 
$\Sigma$ are irreducible and their 
intersections with each of the subspaces $\{ a_j=0\}$ are their proper 
subsets. Therefore the restriction $a_n\neq 0$ and $a_{n-1}\neq 0$ can be 
lifted and one 
concludes that $\alpha _k=2$ 
and $\beta _k=3$ for $k\leq n-1$. 
The number $h\in \mathbb{C}$ is arbitrary and for $n>4$ the set of 
polynomials $P^n$ 
without simple roots is a variety in the space of variables $a$ 
of codimension $\geq 3$. Hence the above reasoning is the proof that for 
$n+1$ the claim of the theorem is true if $k\leq n-1$. 


To perform the induction also for $k=n$ and $k=n+1$ we consider the 
discriminant of the family of polynomials 
$P_*^{n+1}:=a_0x^{n+1}+a_1x^n+\cdots +a_{n+1}$. For its discriminant (denoted 
also by $R^{n+1}$) one has  
$R^{n+1}=(a_0)^{2n}\prod _{1\leq i<j\leq n+1}(z_i-z_j)^2$ ($z_i$ being the roots of 
$P_*^{n+1}$, see \cite{Wi}). Consider the polynomial 
$P^{n+1}_r:=x^{n+1}P_*^{n+1}(1/x)$ (the index $r$ stands for ``reverted''). Its 
roots equal $1/z_i$. Hence its discriminant $R^{n+1}_r$ equals 

$$(a_{n+1})^{2n}\prod _{1\leq i<j\leq n+1}(1/z_i-1/z_j)^2=
(a_0)^{2n}\prod _{1\leq i<j\leq n+1}(z_i-z_j)^2=R^{n+1}~.$$
For $P^{n+1}_r$ the coefficient $a_0$ plays the same role as $a_{n+1}$ plays for 
$P^{n+1}$. Denote by $\tilde{\alpha}_k$, $\tilde{\beta}_k$ the quantities 
$\alpha _k$, $\beta _k$ when defined for the polynomial $P^{n+1}_r$ instead of 
$P^{n+1}$. Hence one can make a shift $x\mapsto x+\tilde{h}$, where $\tilde{h}$ 
is a simple root of $P^{n+1}_r$, and in the same way as above conclude that 
$\tilde{\alpha}_k=2$ and $\tilde{\beta}_k=3$ for $k\leq n-1$. This is 
tantamount to $\alpha _k=2$ and $\beta _k=3$ for $k\geq 2$. As $n\geq 4$, this 
means in particular that $\alpha _n=\alpha _{n+1}=2$ and 
$\beta _n=\beta _{n+1}=3$.

The polynomials $\tilde{D}_k$ and $V_k$ are determinants of Sylvester matrices 
defined after polynomials with integer coefficients. Hence 
$\tilde{D}_k$ and $V_k$ have also integer coefficients. Hence the polynomials 
$M_k$ can also be chosen with integer coefficients 
which implies $c_k\in \mathbb{Q}^*$. 

\end{proof}

\end{document}